\documentclass[12pt]{article}

\usepackage[T1]{fontenc}
\usepackage[breaklinks=true]{hyperref}
\hypersetup{
colorlinks=true,
linkcolor=blue,
urlcolor=blue,
citecolor=black
}
\usepackage{authblk}
\usepackage{ulem}
\usepackage{amsmath,amsfonts,amssymb,lmodern,geometry,enumerate}
\usepackage[font=small,labelfont=bf]{caption}
\usepackage{subcaption}
\usepackage[T1]{fontenc}
\usepackage[english]{babel}
\usepackage{lmodern}
\usepackage{scalefnt}
\usepackage{dsfont}
\usepackage{stmaryrd}
\usepackage{color}
\usepackage{bm,bbm}
\usepackage{mathrsfs,url,color}
\usepackage{breakcites} %for breaking the citations
\usepackage{textcase}
\usepackage{graphicx}
\usepackage{wrapfig}
\usepackage{flushend,cuted}
\usepackage{bm}
\usepackage{tabularx}
\usepackage{color}
\usepackage{indentfirst}
\usepackage{amssymb}
\usepackage{xparse}
\usepackage{pgfplots}

\pgfplotsset{compat=newest}
\usepackage{mdwlist}
\usepackage{amsmath}
\usepackage{amsthm}
\usepackage{nameref}
\makeatletter
\def\@noindentfalse{\global\let\if@noindent\iffalse}
\def\@noindenttrue {\global\let\if@noindent\iftrue}
\def\@aftertheorem{%
  \@noindenttrue
  \everypar{%
    \if@noindent%
      \@noindentfalse\clubpenalty\@M\setbox\z@\lastbox%
    \else%
      \clubpenalty \@clubpenalty\everypar{}%
    \fi}}

\newtheorem{thm}{Theorem}[section]
\AfterEndEnvironment{thm}{\@aftertheorem}

\AfterEndEnvironment{prop}{\@aftertheorem}
\newtheorem{defi}[thm]{Definition}
\AfterEndEnvironment{defi}{\@aftertheorem}
\newtheorem{lma}[thm]{Lemma}
\AfterEndEnvironment{lma}{\@aftertheorem}

\AfterEndEnvironment{cor}{\@aftertheorem}
\newtheorem{re}[thm]{Remark}
\AfterEndEnvironment{re}{\@aftertheorem}
\newtheorem{ap}{Application}

\theoremstyle{remark}
\newcommand{\card}{{\rm card}}

\newcommand{\law}{\mathscr{L}}

\newcommand{\IR}{\mathbb{R}}
\newcommand{\Pro}{\mathbb{P}} 
\newcommand{\prob}{\Pro}
\newcommand{\E}{\mathbb{E}}
\newcommand{\mean}{\E}
\newcommand{\real}{\mathbb{R}}

\newcommand{\Var}{\mathrm{Var}}

\newcommand{\bone}{{\bf 1}}

\newcommand{\cC}{{\cal C}}

\newcommand{\cE}{{\cal E}}

\newcommand{\cN}{{\cal N}}

\def\sS{{\mathscr{S}}}
\def\sX{{\mathscr{X}}}

\def\[{\left[}
\def\]{\right]}
\def\({\left(}
\def\){\right)}

\newcommand{\cov}{{\rm Cov}}
\def\var{{\rm Var}}

\def\ignore#1{}
\newcommand{\Refer}[1]{(\ref{#1})}
\def\qed{\hfill\hbox{${\vcenter{\vbox{
					\hrule height 0.4pt\hbox{\vrule width 0.4pt height 6pt
						\kern5pt\vrule width 0.4pt}\hrule height 0.4pt}}}$}}

\newcounter{con}%for positive constants
\stepcounter{con}
\newcommand{\qcon}[1]{\addtocounter{con}{1}}
\def\sP{{\mathscr{P}}}

\def\given{\typeout{Command 'given' should only be used within bracket command}}
\newcounter{@bracketlevel}
\def\@bracketfactory#1#2#3#4#5#6{
\expandafter\def\csname#1\endcsname##1{%
\addtocounter{@bracketlevel}{1}%
\global\expandafter\let\csname @middummy\alph{@bracketlevel}\endcsname\given%
\global\def\given{\mskip#5\csname#4\endcsname\vert\mskip#6}\csname#4l\endcsname#2##1\csname#4r\endcsname#3%
\global\expandafter\let\expandafter\given\csname @middummy\alph{@bracketlevel}\endcsname
\addtocounter{@bracketlevel}{-1}}%
}
\def\bracketfactory#1#2#3{%
\@bracketfactory{#1}{#2}{#3}{relax}{1mu plus 0.25mu minus 0.25mu}{0.6mu plus 0.15mu minus 0.15mu}
\@bracketfactory{b#1}{#2}{#3}{big}{1mu plus 0.25mu minus 0.25mu}{0.6mu plus 0.15mu minus 0.15mu}
\@bracketfactory{bb#1}{#2}{#3}{Big}{2.4mu plus 0.8mu minus 0.8mu}{1.8mu plus 0.6mu minus 0.6mu}
\@bracketfactory{bbb#1}{#2}{#3}{bigg}{3.2mu plus 1mu minus 1mu}{2.4mu plus 0.75mu minus 0.75mu}
\@bracketfactory{bbbb#1}{#2}{#3}{Bigg}{4mu plus 1mu minus 1mu}{3mu plus 0.75mu minus 0.75mu}
}
\bracketfactory{clc}{\lbrace}{\rbrace}
\bracketfactory{clr}{(}{)}
\bracketfactory{cls}{[}{]}
\bracketfactory{abs}{\lvert}{\rvert}
\bracketfactory{norm}{\Vert}{\Vert}
\bracketfactory{floor}{\lfloor}{\rfloor}
\bracketfactory{ceil}{\lceil}{\rceil}
\bracketfactory{angle}{\langle}{\rangle}

\numberwithin{equation}{section}

\ignore{\makeatletter
	\renewcommand\section{\@startsection {section}{1}{\z@}%
		{-3.5ex \@plus -1ex \@minus -.2ex}%
		{1.3ex \@plus.2ex}%
		{\center\small\sc\MakeTextUppercase}}
	\def\subsection#1{\@startsection {subsection}{2}{0pt}%
		{-3.5ex \@plus -1ex \@minus -.2ex}%
		{1ex \@plus.2ex}%
		{\bf\mathversion{bold}}{#1}}
	\def\subsubsection#1{\@startsection{subsubsection}{3}{0pt}%
		{\medskipamount}%
		{-10pt}%
		{\normalsize\itshape}{\kern-2.2ex. #1.}}
	\makeatother}
\allowdisplaybreaks[4]

\begin{document}
	
\title{\sc\bf\large\MakeUppercase{
			On the rate of normal approximation for Poisson continuum percolation
}}

\author[1]{Tiffany Y.\ Y.\ Lo\thanks{{\sf{email: yin\_yuan.lo@math.uu.se}}. Work supported by the Australian Research Council Grant No DP190100613, Knut and Alice Wallenberg Foundation, Ragnar S\"oderberg Foundation and Swedish Research Council. }}
\author[2]{Aihua Xia\thanks{{\sf{email: aihuaxia@unimelb.edu.au}}. Work supported by the Australian Research Council Grant No DP190100613.}}

%\affil[1]{%
\affil[1]{Department of Mathematics, Uppsala University, L\"agerhyddsv\"agen 1, 752 37 Uppsala, Sweden}
\affil[2]{%
School of Mathematics and Statistics, the University of Melbourne, Parkville VIC 3010, Australia} 

\date{\today}
	
\maketitle
\vskip-1cm
\begin{abstract}
It is known that the number of points in the largest cluster of a percolating Poisson process restricted to a large finite box is asymptotically normal. In this note, we establish a rate of  convergence for the statement. As each point in the largest cluster is determined by points as far as the diameter of the box, known results in the literature of normal approximation for Poisson functionals cannot be directly applied. To disentangle the long-range dependence of the largest cluster, we use the fact that the second largest cluster has comparatively shorter range of dependence to restrict the range of dependence, apply a recently established result in \cite{{CRX20}} to obtain a Berry-Esseen type bound for the normal approximation of the number of points belonging to clusters that have a restricted range of dependence, and then estimate the gap between this quantity and the number of points in the largest cluster.
\end{abstract}

	
	
\section{Introduction and the main result}

Let $\real^m$ be the $m$-dimensional Euclidean space equipped with the Euclidean norm $\|\cdot\|$.
For each $A,B\subset \real^m$, we define $d(A,B)=\inf\{\|x-y\|:\ x\in A,y\in B\}$, where $\inf\emptyset:=\infty$. We write $d(\{x\},B)=:d(x,B)$ for simplicity. Given $r>0$, we define $B(A,r)=\{y\in\real^m:\ d(y,A)<r\}$ and write $B(\{x\},r)=B(x,r)$, so that $B(x,r)$ is simply a ball of radius $r$ with its centre at $x$. We say that a Borel set $A\subset\IR^m$ is {\it connected} with radius $r$ if for any $x_1,x_2\in A$, there exist a finite positive integer $k\ge 2$ and $\{y_1:=x_1,y_2,\ldots,y_{k-1},y_k:=x_2\}\subset A$ such that $B(y_i,r)\cap B(y_{i+1},r)\ne\varnothing$ for all $i=1,\dots,k-1$. We use $\card(A)$ or $|A|$ to denote the cardinality of the set $A$ and for convenience, we use the terms cardinality and size interchangeably. 

\begin{defi} For a Borel set $\sS\subset \real^m$, a subset $A\subset\sS$ is called a {\emph cluster} of $\sS$ with radius $r$ if $A$ is connected with radius $r$ and $B(A,r)\cap B(\sS\setminus A,r)=\varnothing$. 
\end{defi}	

For fixed $r>0$, we say that a point process $\sX$ {\it percolates with radius $r$} if $\sX$ almost surely contains a unique infinite cluster with radius $r$. Let $\sP^\lambda$ be the homogeneous Poisson point process on $\real^m$ with rate $\lambda>0$. For any $r>0$ and $m\geq 2$, it is well known that there exists $0<\lambda_c(r,m)<\infty$ such that $\sP^\lambda$ percolates if and only if $\lambda>\lambda_c(r,m)$; see \cite{ZS85a,ZS85b}. Since $\sP^\lambda$ is scale invariant, it is enough to consider only $r=1$ and write $\lambda_c:=\lambda_c(m):=\lambda_c(1,m)$. From now on, any cluster with radius 1 is simply referred to as a cluster. Proving the exact values of $\lambda_c(m)$ remains an open question, although for $m=2$, a sharp estimate was given in \cite{BBW05}. There is also a vast literature considering more general continuum percolation since \cite{G61} initiated the study, where the point processes are not necessarily homogeneous Poisson and each ball $B(x,r)$ can be replaced by a random shape centred at~$x$; we refer the reader to \cite{MR96} for a comprehensive overview.

In practice, any physical system is finite and the percolation phenomenon is examined through growing observation windows, hence it is of practical interest to study the asymptotic behaviour of the cardinality of the largest cluster inside a growing window in $\mathbb{R}^m$. The statistical behaviour of the largest cluster in a large finite observation window under both the regimes $\lambda<\lambda_c$ and $\lambda>\lambda_c$ have been thoroughly investigated in \cite{Pe03,PP96}, and here we briefly summarise some results for the case $\lambda>\lambda_c$. Let $\Gamma_n:=[-n/2,n/2]^m$ be such a window and $N_n$ be the number of points in the largest cluster in $\sP^\lambda_n:=\sP^\lambda\cap \Gamma_n$. It was shown in \cite[Chapter 10]{Pe03} that when $m\geq 2$, $n^{-m} N_n\to \lambda p(\lambda)$ in probability as $n\to\infty$, where $0<p(\lambda)<1$ is the probability that the infinite cluster contains the origin. Furthermore, with probability tending to one as $n\to\infty$, the size of the second largest cluster in $\sP^\lambda_n$ is of the exact order $\Theta\left((\ln n)^{m/(m-1)}\right)$, thus establishing the uniqueness of the largest cluster. Large deviation estimates for the size, volume, and diameter of the largest cluster were provided in \cite{PP96}. The result that is most pertinent to our work here is the central limit theorem for $N_n$  established in \cite[Theorem 10.22]{Pe03} that holds for $m\geq 2$ and $\lambda>\lambda_c$. Let $B_n^2=\Var(N_n)$. \cite[Theorem~10.22]{Pe03} and the errata in \cite{Pe10} showed that, for $m\geq 2$, there exists a constant $0<\sigma^2:=\sigma^2(\lambda,m)<\infty$ such that  
\begin{equation}\label{eq:variance}
    n^{-m} B_n^2\to \sigma^2
\end{equation}
and 
\begin{equation*}
   \law( n^{-m/2}(N_n-\E N_n)) \to \cN(0,\sigma^2)\quad \text{as $n\to\infty$.}
\end{equation*}

Our main result below compliments this central limit theorem by providing a convergence rate in the Kolmogorov distance.
%{\color{blue} This variance is for the size of the largest component, not the component carved out from the infinite component, which can include the smaller islands in the $\Gamma$. The theorem below assumes that the variance is of the same order as above.}

\begin{thm}\label{mainthm1}
Suppose that $m\geq 2$ and $\lambda>\lambda_c$, and let $W_n:=(N_n-\E N_n)/B_n$. Then 
\begin{equation}\label{mainthm1.1}
    d_{\mathrm{K}}(\law(W_n),\cN(0,1)):=\sup_{x\in\real}|\prob(W_n\le x)-\prob(Z\le x)|     \le O\left(n^{-m/2}(\ln n)^{2m}\right),
\end{equation}
where $Z\sim \cN(0,1)$.
\end{thm}

\begin{re} The logarithmic factor in \Refer{mainthm1.1} seems unavoidable because the percolation is a long-range dependent structure and it differs significantly from the geometric structures studied in \cite{SY21}. However, we suspect that the dependence of dimensionality in $(\ln n)^{2m}$ is due to the choice of the window size that we use to construct another score function with local dependence and it is not clear whether one can reduce or remove the dependence on $m$ with a smaller window size or another method.
\end{re}

The proof of the central limit theorem in \cite{Pe03} hinges on a martingale argument, while here we rely on Stein's method \cite{CGS11} to deduce the convergence and the rate in Theorem~\ref{mainthm1}. 
Besides Stein method, one may also consider other tools such as the stabilisation tool \cite{PY01, PY05}, the Malliavin-Stein technique via the Wiener-It\^{o} expansion \cite{PSTU10} and the second order Poincar\'e inequalities \cite{LPS16}. In fact, using these tools, a variety of central limit theorems have been developed for random quantities of the form $\sum_{x\in \sP^\lambda \cap A}\xi(x,\sP^\lambda)$, where $A\subset \mathbb{R}^m$ is a bounded Borel set, and $\xi(x,\sP^\lambda)$ is a score function that measures the contribution of $x$ with respect to $\sP^\lambda$; see for examples \cite{BX, PY01, PY03, PY05, PSTU10, S12, S16, LPS16, LSY19, CX20, SY21, BM22}. In our setting, $A$ is replaced by $\Gamma_n$, and $\xi(x,\sP^\lambda)$ can be taken as the indicator function that takes value one if the Poisson point $x$ belongs to the largest cluster in $\sP^\lambda_n$. To obtain rates of convergence, the existing literature on normal approximation generally requires the score functions to have short-range dependence, which loosely speaking, is the condition that the score function $\xi(x,\sP^\lambda)$ depends only on points that are not too far away from $x$. For instance, \cite{BM22, Pe07, S16, XY15, SY21} require the score function $\xi(x,\sP^\lambda)$ to be determined by the points of $\sP^\lambda$ in a region near $x$ or a ball $B(x,R)$
with a random radius $R$ such that $\prob(R>t)$ decreases as the reciprocal of a polynomial or an exponential function of $t$ as $t\to\infty$. In our case, with the long-range dependence of the points in the percolation, $\prob(R=\Theta(n))\approx1$, so the score function in consideration here does not fit into the framework of such literature.

\vspace{0.5cm}
{\noindent\textit{Strategy of the proof.} To disentangle the long-range dependence, we use the characteristic of the second largest cluster to construct a suitable score function $\xi'(x,\sP^\lambda)$ that takes value one if $x$ belongs to a `local' cluster that is typically larger than the second largest cluster, apply \cite[Corollary~3.2]{CRX20} to obtain a Berry-Esseen type bound for the normal approximation of the sum $N_{\theta,n}'$ of these score functions, and then bound the gap between $N_n$ and $N_{\theta,n}'$. }
 
\section{The proof of Theorem~\ref{mainthm1}}

To represent $N_n$ as the sum of appropriate score functions, for any $\sX\subset\real^m$, we 
write $\cC(\sX)$ as the set of all clusters of $\sX$, and for $x\in\sX$, let {$\cC(x,\sX)$} be the cluster of $\sX$ containing $x$, and write $\cC_0(\sX)$ as the largest cluster of $\sX$ if it is unique. Furthermore, 
define the score function of the point configuration $\sP^\lambda_n$ at $x$ as $\xi(x,\sP^\lambda_n):={\bf 1}_{[\cC(x,\sP^\lambda_n)=\cC_0(\sP^\lambda_n)]}$. The score function collects the points in the largest cluster in $\sP^\lambda_n$ and $N_n=\sum_{x\in\sP^\lambda_n}\xi(x,\sP^\lambda_n)$.  
%Additionally, define $\Xi(dx):=\xi(x,\sP^\lambda_n)\sP^\lambda_n(dx)$.

To tackle the long-range dependence, we first observe that the typical size of the second largest cluster in {$\sP^\lambda_n$} is no more than $c(\ln n)^{m/(m-1)}$ for a constant $c>0$ not depending on $n$. Next, for each $x\in \real^m$, we take the cube $A_{x,\theta,1}$ with the centre $x$ and edge length $2\theta\ln n$, {$A_{x,\theta,2}=A_{x,\theta,1}\cap \Gamma_n$}, and show that the point $x$ is in the largest cluster in {$\sP^\lambda_n$} is essentially the same as that $\cC(x,\sP^\lambda\cap A_{x,\theta,2})$ is the largest cluster in $\sP^\lambda\cap A_{x,\theta,2}$. However, the latter characterisation ensures that its corresponding score function has short-range dependence, so the tools of normal approximation to the sum of locally dependent score functions can be applied.  For the size of the second largest cluster, the following lemma is a direct consequence of modifying (10.56) and (10.58) in the proof of \cite[Theorem 10.18]{Pe03}; noting that the proof itself is an application of \cite[Theorem~2]{PP96}.

\begin{lma}
Let $\cC^{(2)}_0(\sP^\lambda_n)$ be the second largest cluster in $\sP^\lambda_n$. Then there exists $k_0>0$ such that for any $k_1\ge k_0$, there are $k_2(k_1)=:k_2>0$ and $n_0(k_1)=:n_0>0$ not depending on $n$ such that 
\begin{equation}\label{eq:2ndcl}
\prob\left(|\cC^{(2)}_0(\sP^\lambda_n)|\geq k_1 (\ln n)^{m/(m-1)}\right)\leq k_2n^{-10m}
\end{equation}
for $n\ge n_0$.
\end{lma}

Write $\cC_0(\sP^\lambda\cap A_{x,\theta,2})$ as the largest cluster in $\sP^\lambda\cap A_{x,\theta,2}$ if it is unique, and define another score function $\xi'(x,\theta,\sP^\lambda_n):={\bf 1}_{[\cC(x,\sP^\lambda\cap A_{x,\theta,2})=\cC_0(\sP^\lambda\cap A_{x,\theta,2})]}$, so that $\xi'(x,\theta,\sP^\lambda_n)=1$ if $x$ belongs to the largest cluster of $\sP^\lambda\cap A_{x,\theta,2}$. We now assess the difference between $N_n$ and $N_{\theta,n}':=\sum_{x\in\sP^\lambda_n}\xi'(x,\theta,\sP^\lambda_n)$.

\begin{lma}\label{lma2} There exists a constant $\theta>0$ such that
\begin{align*}
    \prob\left(N_{\theta,n}'-N_n\ne0\right)&=O\left(n^{-10m}\right).
\end{align*}
\end{lma}
\noindent {\it Proof.} Let $\tilde\theta:=\tilde\theta(\lambda)$ denote the probability that there is an unbounded cluster $D$ such that $B(D,1)$ intersects the {ball} of unit volume centred at the origin ${\bf 0}\in \real^m$. Furthermore, let $\cE_{1}$ be the event that the largest cluster $\cC_0(\sP^\lambda_n)$ is the unique cluster such that $|\cC_0(\sP^\lambda_n)|\geq0.5\lambda\tilde\theta n^m$ and with diameter at least $0.5n$, where the diameter of a {subset} $A\subset \IR^m$ is $\sup\{\|x-y\|: x,y\in A\}$. Then, \cite[Theorem~2]{PP96} states that there exist constants $k_3>0$ and $n_1>0$ such that
\begin{align}
\prob\left(\cE_1^c\right)&\le \exp\{-k_3 n\},\ \ \ n\ge n_1.\label{diff01}
\end{align}
For any {$x\in \Gamma_n$}, let $\cE_{1,x}$ be the counterpart of $\cE_{1}$ with $\sP^\lambda_n$ replaced with $\sP^\lambda_n\cup \{x\}$.  Since the extra point $x$ does not reduce the largest cluster, \Refer{diff01} implies that
\begin{equation}
\prob\left(\cE_{1,x}^c \right)\le \exp\{-k_3 n\},\ \ \ n\ge n_1.\label{diff02}
\end{equation}

 Let $k_1$, $k_2$ and $k_3$ be as in \Refer{eq:2ndcl} and \Refer{diff01}, $\cE_2:=\{|\cC^{(2)}_0(\sP^\lambda_n)|< k_1 (\ln n)^{m/(m-1)}\}$ and $\theta=11m/ k_3$. In addition, for any $x\in{\sP^\lambda_n}$, let $E_{0,x}$ be the event that the largest cluster in $\sP^\lambda\cap A_{x,\theta,2}$ is unique, $|\cC_0(\sP^\lambda\cap A_{x,\theta,2})|\geq 0.5\lambda\tilde \theta (\theta\ln n)^m$ and it is of diameter at least {$0.5\theta \ln n$}. We claim that
\begin{equation}\label{eq:goodevents}
    E_{0,x}\cap \cE_1 \cap \cE_2 \subset \clc{\xi(x,\sP^\lambda_n)=\xi'(x,\theta,\sP^\lambda_n)},
\end{equation}
or equivalently $\{\xi(x,\sP^\lambda_n)\ne\xi'(x,\theta,\sP^\lambda_n)\}\cap E_{0,x}\cap \cE_1 \cap \cE_2 =\emptyset$.
 
We first consider the case where $x \in \sP^\lambda_n$ belongs to $\cC_0(\sP^\lambda_n)$ but not $\cC_0(\sP^\lambda_n\cap A_{x,\theta,2})$, i.e.\ $\{\xi(x,\sP^\lambda_n)=1,\xi'(x,\theta,\sP^\lambda_n)=0\}$ or equivalently, $\{\cC(x, \sP^\lambda_n)=\cC_0(\sP^\lambda_n),\cC(x,\sP^\lambda\cap A_{x,\theta,2})\ne\cC_0(\sP^\lambda\cap A_{x,\theta,2})\}$. Since $x$ is the centre of {$A_{x,\theta,1}$} and $x$ belongs to $\cC_0(\sP^\lambda_n)$, $\cC_0(\sP^\lambda_n)\cap A_{x,\theta,2}$ contains $x$ and $\cC(x, \sP^\lambda\cap A_{x,\theta,2})$ must have diameter at least $\theta\ln n-1$. On the event $E_{0,x}\cap \cE_1$, the only cluster with diameter at least $0.5\theta \ln n$ is $\cC_0(\sP^\lambda\cap A_{x,\theta,2})$, and so $\cC(x,\sP^\lambda\cap A_{x,\theta,2})=\cC_0(\sP^\lambda\cap A_{x,\theta,2})$ and $\xi'(x,\theta,\sP^\lambda_n)=1$, which is in contradiction to $\xi'(x,\theta,\sP^\lambda_n)=0$.
%On the other hand, if $x$ does not belong to $\cC_0(\sP^\lambda\cap A_{x,\theta,2})$, then $\cC(x,\sP^\lambda\cap A_{x,\theta,2})$ must have diameter less than $\theta\ln n-1$. Thus, $\cC(x,\sP^\lambda\cap A_{x,\theta,2})$ cannot intersect $\partial A_{x,\theta,2} \oplus B(1)$, with B(1) being the unit $m$-dimensional ball and so it cannot intersect the discs of the points belonging to $\cC_0(\sP^\lambda)\cap (\Gamma\setminus A_{x,\theta,2})$. Consequently, $\{\xi(x,\sP^\lambda_n)=1,\xi'(x,\theta,\sP^\lambda_n)=0\}$ cannot be possible.

We turn to the other case where $x \in \sP^\lambda_n$ belongs to $\cC_0(\sP^\lambda_n\cap A_{x,\theta,2})$ but not $\cC_0(\sP^\lambda_n)$, i.e.\ $\{\xi(x,\sP^\lambda_n)=0,\xi'(x,\theta,\sP^\lambda_n)=1\}$. On the event $E_{0,x}\cap\cE_1\cap \cE_2$, the second largest cluster $\cC^{(2)}_0(\sP^\lambda_n)$ has at most $k_1(\ln n)^{m/(m-1)}$ points, while the largest cluster in $A_{x,\theta,2}\cap \sP^\lambda$ has at least $0.5\lambda\tilde \theta (\theta\ln n)^m$ points, hence if $\xi'(x,\theta,\sP^\lambda_n)=1$ and $\xi(x,\sP^\lambda_n)=0$, then $\cC(x,\sP^\lambda\cap A_{x,\theta,2})=\cC_0(\sP^\lambda\cap A_{x,\theta,2})$ is no larger than $\cC^{(2)}_0(\sP^\lambda_n)$, giving  $|\cC_0(\sP^\lambda\cap A_{x,\theta,2})|\leq k_1(\ln n)^{m/(m-1)}${, which leads to a contradiction.} This concludes the proof of \Refer{eq:goodevents}.

{Let $\cE_0:=\cap_{x\in \sP^\lambda_n}E_{0,x}$}. We have $\cap_{x\in\sP^\lambda_n} \clc{\xi(x,\sP^\lambda_n)=\xi'(x,\theta,\sP^\lambda_n)} \subset \{N_n= N'_{\theta,n}\}$, {and} by \Refer{eq:goodevents},
\begin{align}\label{eq:finalbd}
    \prob(N_n\ne N'_{\theta,n})\leq \prob(\cup_{x\in\sP^\lambda_n}\clc{\xi(x,\sP^\lambda_n)\ne \xi'(x,\theta,\sP^\lambda_n)})\leq \prob(\cE_0^c) + \prob(\cE^c_1) + \prob(\cE^c_2).
\end{align}
By \Refer{diff01} and Lemma~\ref{eq:2ndcl}, $\prob(\cE_1^c)\leq \exp\{-k_3n\}$ and $\prob(\cE_2^c)\leq k_2n^{-10m}$. Using the Palm distributions of Poisson $\sP^\lambda$ \cite[Chapter~10]{Kallenberg83} and \Refer{diff02}, for $n\ge e^{n_1k_3/(11m)}$,
\begin{align*}
    \prob(\cE_0^c)\leq \E\int_{{\Gamma_n}} {\bf 1}_{E^c_{0,x}} \sP^\lambda (dx) = \lambda \int_{{\Gamma_n}} \prob(E^c_{0,x}) dx \leq O(n^m \exp\{-k_3\theta \ln n\})= O(n^{-10m}),
\end{align*}
where the last equality follows from $\theta=11m/k_3$. Applying these bounds to \Refer{eq:finalbd} concludes the proof.
\qed

\begin{lma}\label{lma3} With $\theta$ as in the proof of Lemma~\ref{lma2}, we have
\begin{align*}
    \mean|N'_{\theta,n}-N_n|&\le O\left(n^{-4m}\right),\\
    |\var(N_{\theta,n}')-\var(N_n)|&=O\left(n^{-m}\right),\\
    \var(N_{\theta,n}')&=\Theta(n^m).
\end{align*}
\end{lma}

\noindent {\it Proof.} For the first claim, {below} we apply the Cauchy-Schwarz inequality in the second inequality and Lemma~\ref{lma2} in the last inequality to get
\begin{align*}
\mean|N'_{\theta,n}-N_n|&\le \mean\left(|\sP^\lambda_n|\bone_{[N_{\theta,n}'\ne N_n]}\right)\\
&\le \sqrt{\mean(|\sP^\lambda_n|^2)\prob(N_{\theta,n}'\ne N_n)}\le O\left(n^m\right)O\left(n^{-5m}\right)=O\left(n^{-4m}\right).
\end{align*}

Likewise, we have 
\begin{align*}\var(N'_{\theta,n}-N_n)&\le \mean[(N'_{\theta,n}-N_n)^2]\le \mean\left(|\sP^\lambda_n|^2\bone_{[N_{\theta,n}'\ne N_n]}\right)\\
&\le \sqrt{\mean(|\sP^\lambda_n|^4)\prob(N_{\theta,n}'\ne N_n)}\le O\left(n^{2m}\right)O\left(n^{-5m}\right)=O\left(n^{-3m}\right),
\end{align*}
therefore,
\begin{align*}|\var(N'_{\theta,n})-\var(N_n)|&= |\var(N'_{\theta,n}-N_n)+2\cov(N_n,N_{\theta,n}'-N_n)|\\
&\le O\left(n^{-3m}\right)+2\sqrt{\var(N_n)\var(N_{\theta,n}'-N_n)}=O\left(n^{-m}\right).
\end{align*}

The third claim follows from \Refer{eq:variance} and the second claim. It can also be directly obtained from \cite[Lemma~4.6]{XY15} and the fact that the score function $\xi'$ is locally dependent.  \qed

We now establish the error bound of the normal approximation to $N'_{\theta,n}$. Let $B'_{\theta,n}$ be the standard deviation of $N'_{\theta,n}$.

\begin{lma}\label{lma4} For any constant $\theta>0$, let {$W_{\theta,n}'=(N'_{\theta,n}-\mean N'_{\theta,n})/B'_{\theta,n}$}, then 
$$
    d_{\mathrm{K}}(\law(W_{\theta,n}'),\cN(0,1))= O\left( (\ln n)^{2m}n^{-m/2}\right).
$$
\end{lma}

\noindent {\it Proof.}  Recall that $\xi'(x,\theta,\sP^\lambda_n)={\bf 1}_{[\cC(x,\sP^\lambda\cap A_{x,\theta,2})=\cC_0(\sP^\lambda\cap A_{x,\theta,2})]}$, define the point process $\Xi'(dx)=\xi'(x,\theta,\sP^\lambda_n)\sP^\lambda_n(dx)$, $\Lambda'(dx)=\mean \Xi'(dx)$, let $\Xi'_x$ be its Palm process at~$x$ \cite[Chapter~10]{Kallenberg83}. Let $A_{x,\theta,3}$ be the cube with centre $x$ and edge length $4\theta\ln n$ and $A_{x,\theta,4}:=A_{x,\theta,3}\cap \Gamma_n$. Because the score function $\xi'(x,\theta,\sP^\lambda_n)$ is completely determined by the point configuration $\sP^\lambda\cap A_{x,\theta,2}$, we can construct $\Xi'$ and $\Xi'_x$ together such that $\Xi'_x$ and $\Xi'$ are identical outside $A_{x,\theta,4}$. Let $Y'_x=\Xi'_x(\Gamma_n)-\Xi'(\Gamma_n)=\Xi'_x(A_{x,\theta,4})-\Xi'(A_{x,\theta,4})$, $\Delta'_x=Y'_x/B'_{\theta,n}$, $D=\{(x,y)\in\Gamma_n^2:\ d(x,y)\le 4\sqrt{m}\theta\ln n\}$, then
$\Delta'_x$ and $\Delta'_y$ are independent when $(x,y)\in\Gamma_n^2\setminus D$. By \cite[Corollary~3.2]{CRX20},
$$
d_{\mathrm{K}}(\law(W'_{\theta,n}),\cN(0,1)) \le 7s_1 + 5.5s_2 + 10s_3,
$$
where, with $B:=B'_{\theta,n}$,
\begin{align*}
s_1 &= \frac{1}{B^2}\left(\int_{(x,y)\in D}\E\{(Y'_x)^2 \bone_{[|Y'_x|\le B]}\}\Lambda'(dx)\Lambda'(dy)\right)^{\frac{1}{2}};\\
s_2 &= \frac{1}{B^3}\int_{\Gamma_n}\E \{(Y'_x)^2\}\Lambda'(dx); \\
s_3 &= \frac{1}{B^3}\int_{(x,y)\in D} \E \{|Y'_x|\bone_{[|Y'_x|\le B]}\}\Lambda'(dx)\Lambda'(dy).
\end{align*}
Since $|Y'_x|\le \sP^\lambda(A_{x,\theta,4})$, we have $\E |Y'_x|\le\lambda (4\theta\ln n)^m$ and 
$$\E \{(Y'_x)^2\}\le {2\cdot \lambda^2 (4\theta\ln n)^{2m}}$$ for large $n$. Hence,
\begin{align*}
s_1 &\le \frac{1}{B^2}\left(\int_{(x,y)\in D}{2\lambda^2 (4\theta\ln n)^{2m}}\Lambda'(dx)\Lambda'(dy)\right)^{\frac{1}{2}}\le O\left(\frac{n^{m/2}(\ln n)^{1.5m}}{B^2}\right);\\
s_2 &\le \frac{1}{B^3}\int_{\Gamma_n} {2\lambda^2(4\theta\ln n)^{2m}}\Lambda'(dx)=O\left(\frac{n^m(\ln n)^{2m}}{B^3}\right);\\
s_3 &\le \frac{1}{B^3}\int_{(x,y)\in D} {\lambda (4\theta\ln n)^{m}}\Lambda'(dx)\Lambda'(dy)=O\left(\frac{n^m(\ln n)^{2m}}{B^3}\right).
\end{align*}
By Lemma~\ref{lma3}, we have $B^2=\Theta(n^m)$, the proof is complete.
\qed

With these preparations, we are now ready to prove Theorem~\ref{mainthm1}.

\noindent {\it Proof of Theorem~\ref{mainthm1}.} Using the triangle inequality, we have
$$d_{\mathrm{K}}\left(\law(W_n),\cN(0,1)\right)\le d_{\mathrm{K}}\left(\law(W_n),\law(W_{\theta,n}')\right)+d_{\mathrm{K}}\left(\law(W_{\theta,n}'),\cN(0,1)\right),$$
hence, by Lemma~\ref{lma4}, it suffices to show that
\begin{equation}d_{\mathrm{K}}\left(\law(W_n),\law(W_{\theta,n}')\right)\le 2d_{\mathrm{K}}\left(\law(W_{\theta,n}'),\cN(0,1)\right)+O\left(n^{-2m}\right).\label{proofofmainthm01}\end{equation}
To this end, let {$V_{\theta,n}:=(N'_{\theta,n}-\mean N_n)/B_n$}, $v_{\theta,n}=\mean V_{\theta,n}$, $r_{\theta,n}^2:=\var(V_{\theta,n})=(B_{\theta,n}/B_n)^2$, applying the triangle inequality in the first inequality and \cite[(5.9)]{XY15} in the second inequality below, we have
\begin{align}
&d_{\mathrm{K}}(\law(W_n),\law(W_{\theta,n}'))\nonumber\\
&\le d_{\mathrm{K}}\left(\law(W_n),\law(V_{\theta,n})\right)+d_{\mathrm{K}}\left(\law(V_{\theta,n}), \cN\left(v_{\theta,n},r_{\theta,n}^2\right)\right)\nonumber\\
&\ \ \ +d_{\mathrm{K}}\left(\cN\left(v_{\theta,n},r_{\theta,n}^2\right),\cN(0,1)\right)+d_{\mathrm{K}}\left(\law(W_{\theta,n}'),\cN(0,1)\right)\nonumber\\
&\le \prob(N_n\ne N'_{\theta,n})+2d_{\mathrm{K}}\left(\law(W_{\theta,n}'),\cN(0,1)\right)+\frac{|v_{\theta,n}|}{\sqrt{2\pi}}
+\frac{|r_{\theta,n}^2-1|}{\sqrt{2e\pi}}.\label{proofofmainthm02}
\end{align}
Lemma~\ref{lma3} gives
\begin{equation}|v_{\theta,n}|\le O\left(n^{-4.5m}\right),\ \ \
|r_{\theta,n}^2-1|\le O\left(n^{-2m}\right),
\label{proofofmainthm03}\end{equation}
hence \Refer{proofofmainthm01} follows from combining Lemma~\ref{lma2} and the estimates \Refer{proofofmainthm02} and \Refer{proofofmainthm03}.  \qed

\newlength{\bibsep}%{0.5ex}
\def\bibfont{\small}
\normalem

\def\ac{{Academic Press}~}
\def\aap{{Adv. Appl. Prob.}~}
\def\ap{{Ann. Probab.}~}
\def\anap{{Ann. Appl. Probab.}~}
\def\eljp{{\it Electron.\ J.~Probab.\/}~} 
\def\jap{{J. Appl. Probab.}~}
\def\jrss{{J. R. Stat. Soc.}~}
\def\jtp{{J. Theor. Probab.}~}
\def\jws{{John Wiley $\&$ Sons}~}
\def\ny{{New York}~}
\def\ptrf{{Probab. Theory Related Fields}~}
\def\sp{{Springer}~}
\def\spa{{Stochastic Process. Appl.}~}
\def\spl{{Stat. Probab. Lett.}~}
\def\sv{{Springer-Verlag}~}
\def\tpa{{Theory Probab. Appl.}~}
\def\zw{{Z. Wahrsch. Verw. Gebiete}~}

\end{document}